\documentclass[12pt,a4paper]{article}
\usepackage{amssymb}
\usepackage{amsmath}
\usepackage{amsthm}
\usepackage{times} 
\newcommand{\nc}{\newcommand} 
\nc{\at}{\begin{ttt}}\nc{\zt}{\end{ttt}}
\nc{\mb}{\mathbb}
\nc{\A}{\mb{A}}
\nc{\bb}{\bigskip}
\nc{\C}{\mb{C}}
\nc{\cl}{\centerline}
\nc{\ind}{\hskip 1em\relax}

\newtheorem*{ttt}{Theorem}

\begin{document}
\begin{center}
{\Huge The functional equation of the zeta \bb

function of a global field}\bb
\end{center}

All the references are made to\bb

\cl{[W] Andr\'e Weil, {\bf Basic Number Theory}, Springer, 1967.}\bb

\ind Recall that a {\it global field} is a finite degree extension of $\mb{Q}$ or
of $\mb{F}_p(T)$. 
Let $k$ be a global field. The completion $k_v$ of $k$ at a finite
place $v$ has a largest compact subring; this subring is open and local,
and its residual field has finite cardinality $q_v$. The product of the
$(1-q_v^{-s})^{-1}$ converges for $\Re s>1$ to a holomorphic function which
can be continued to a meromorphic function $\zeta_k$ on the whole complex
plane (see [W], Equation (6) p. 128 and Definition~8 p. 129). Put 

$$Z_k(s):=\Big(\pi^{-s/2}\ \Gamma(s/2)\Big)^{r_1}
\Big((2\pi)^{1-s}\ \Gamma(s)\Big)^{r_2}\zeta_k(s),$$

where $r_1$ and $r_2$ denote respectively the numbers of real and imaginary
places of $k$. For each place $v$ of $k$ define the Haar measure $\beta_v$
on $k_v$ as follows. If $k_v$ is nonarchimedian, let $\beta_v$ be the Haar
measure giving the volume 1 to the largest compact subring; if
$k_v=\mb{R}$, put $\beta_v:=dx$; if $k_v=\C$, set 
$\beta_v:=|dz\wedge d\overline{z}|$. The product of the $\beta_v$ defines a
Haar measure $\beta$ on the locally compact $k$-algebra $\A$ of adeles of
$k$. Recall that $k$ is discrete and cocompact in $\A$; in particular the
$\beta$-covolume $\beta(\A/k)$ of $k$ in $\A$ is a positive real number. 

\at We have\bb

\cl{\boxed{Z_k(1-s)=\beta(\A/k)^{2s-1}Z_k(s)}}\zt

{\bf Observations.} If $k$ is a number field, then $\beta(\A/k)$ is the
square root of the absolute value of the discriminant (see Proposition~7
p. 91 of [W]). If $k$ is of positive characteristic, if $q$ is the cardinality
of its largest finite subfield (which exists by Theorem~8 p. 77 of [W]), and
if $g$ is its genus, then we have $\beta(\A/k)=q^{g-1}$  (see
Corollary~1 p. 100 of [W]). \bb

{\bf Proof of the Theorem.} Let $\Phi:\A\to\C$ be the function defined
on page 127 of [W], and let $\A^\times$ be the idele group of $k$. For any
continuous morphism $\omega$ from $\A^\times$ to $\C^\times$ put 

$$Z(\omega,\Phi):=\int_{\A^\times}\Phi(x)\ \omega(x)\ d\mu(x),$$

where $\mu$ is the measure defined on page 120 of [W]. For any $s\in\C$ and
any $x\in\A^\times$ let $\omega_s(x)$ be the adelic norm of $x$ raised to
the power $s$. Since there is a positive constant $c_k$ such that 

$$Z_k(s)=c_k\ Z(\omega_s,\Phi)$$

by Equation~(6) page 128 of [W], it suffices to check

$$Z(\omega_{1-s},\Phi)=\beta(\A/k)^{2s-1}Z(\omega_s,\Phi).$$

Identify $\A$ to its topological dual via a nontrivial character of $\A$
trivial on $k$ as explained in Theorem~3 p. 66 of [W], and equip $\A$ with
its selfdual measure, so that the Fourier transform of a continuous
integrable function on $\A$ is a well defined function on $\A$. Letting
$\Phi'$ be the Fourier transform of $\Phi$, we claim 

$$Z(\omega_{1-s},\Phi)=Z(\omega_s,\Phi')
=\beta(\A/k)^{2s-1}Z(\omega_s,\Phi).$$

Indeed the first equality follows from Theorem~2 p. 121 of [W], and the
second one from the above Observations, Proposition~6 p. 113 of [W], and
Equation~(7) page 129 of [W]. QED\bb\bb

\hfill\small Pierre-Yves Gaillard, Universit\'e Nancy 1, 02.28.05

%\vfill\hfill\tiny050228.zeta, Mon Feb 28 14:43:58 CET 2005

\end{document}